\documentclass[12pt]{amsart}
\usepackage{mathrsfs}
\usepackage{amsmath,amssymb}
\usepackage{amsfonts}
\usepackage{amsthm}
\usepackage{latexsym}
\usepackage{graphicx}
\def\p{\partial}

\def\Hom{\mbox{Hom}}

\def\R{\mathbb{R}}

\def\vv<#1>{\langle#1\rangle}

\def\vv<#1>{\left\langle#1\right\rangle}
\def\pt{\frac{\p}{\p t}}
\def\pu{\frac{\p}{\p u}}

\def\wt{\widetilde}
\def\e{\epsilon}
\newtheorem{thm}{Theorem}[section]

\theoremstyle{definition}
\newtheorem{defn}{Definition}[section]
\theoremstyle{remark}

\newtheorem{rem}{Remark}[section]

\numberwithin{equation}{section}
\def\dev{{\rm dev}}

\begin{document}

\title{Developments of curves with respect to symmetric tensors and existence of isometric immersions with prescribed second fundamental form}

\author{Chengjie Yu$^1$}
\address{Department of Mathematics, Shantou University, Shantou, Guangdong, 515063, China}
\email{cjyu@stu.edu.cn}
\thanks{$^1$Research partially supported by GDNSF with contract no. 2021A1515010264 and NNSF of China with contract no. 11571215.}
\renewcommand{\subjclassname}{%
  \textup{2010} Mathematics Subject Classification}
\date{}
\begin{abstract}
In this paper, we introduce the notion of developments of curves with respect to  symmetric tensors and use it to prove the existence of isometric immersions into a general ambient space with prescribed second fundamental form. Our method  provides a geometric construction of such an isometric immersion.
\end{abstract}
\maketitle\markboth{Chengjie Yu }{Isometric immersions with prescribed second fundamental form}
\section{Introduction}
Let $(M^n,g)$ be a Riemannian submanifold of $(\widetilde M^{n+s},\wt g)$. Then, the curvature tensors $R$ and $\wt R$ for $M$ and $\wt M$ respectively are related by the following three equations:
\begin{equation}\label{eq-Gauss}
R(X,Y,Z,W)=\wt R(X,Y,Z,W)+\vv<h(X,W),h(Y,Z)>-\vv<h(X,Z),h(Y,W)>,
\end{equation}
\begin{equation}\label{eq-Codazzi}
\vv<(\nabla^\perp_{X}h)(Y,Z)-(\nabla^\perp_Yh)(X,Z),\xi>=\wt R(Z,\xi,X,Y),
\end{equation}
and
\begin{equation}\label{eq-Ricci}
R^\perp(\xi,\eta,X,Y)=\wt R(\xi,\eta,X,Y)+\vv<A_\xi(Y),A_\eta(X)>-\vv<A_\eta(Y),A_\xi(X)>
\end{equation}
for any $X,Y,Z,W\in \Gamma(TM)$ and $\xi,\eta\in \Gamma(T^\perp M)$, where
\begin{equation*}
h(X,Y)=\left(\wt\nabla_XY\right)^\perp\  (\forall X,Y\in \Gamma(TM))
\end{equation*}
is the second fundamental form with $\wt\nabla$ the Levi-Civita connection for $\wt g$, $\nabla^\perp$ is the normal connection, $R^\perp$ is the curvature tensor for the normal connection, and $A_\xi(X)$ is the shape operator defined by
\begin{equation*}
\vv<A_\xi(X),Y>=\vv<h(X,Y),\xi>\ (\forall X,Y\in \Gamma(TM),\ \xi\in \Gamma(T^\perp M)).
\end{equation*}
The equations \eqref{eq-Gauss}, \eqref{eq-Codazzi} and \eqref{eq-Ricci} are called the Gauss equation, Codazzi equation and Ricci equation respectively. The classical fundamental theorem for submanifolds tracing back to the work of Bonnet \cite{Bo} says that the converse of the above is true when the ambient space $(\wt M,\wt g)$ is a space form.  For the precise statement of the result, see \cite{Chen}. The reason one requires that the ambient space $(\wt M,\wt g)$ is a space form is that the curvature tensor $\wt R$ is intrinsic in this case. More precisely, when $\wt M$ is a space form of sectional curvature $K$, one has
\begin{equation*}
\wt R(X,Y,Z,W)=K\left(\vv<X,W>\vv<Y,Z>-\vv<X,Z>\vv<Y,W>\right)
\end{equation*}
for any tangent vectors $X,Y,Z,W$ of $\wt M$. There are some efforts to extend the fundamental theorem for submanifolds to more general ambient spaces such as products of space forms, symmetric spaces in the past decades.  See for examples \cite{CC,Da,Ko,LZ,LT} and references there in.

Note that the existence part of the fundamental theorem for submanifolds can be viewed as the existence of isometric immersions into space forms with prescribed second fundamental form. In this paper, we consider the existence of isometric immersion into general ambient spaces with prescribed second fundamental form. Our idea is to view such an existence result as a postive codimensional Cartan-Ambrose-Hicks theorem (See \cite{Am,Ca,CE,Hi}) and use a similar argument of the proof for the Cartan-Ambrose-Hicks theorem to show the existence of such kinds of isometric immersions. Similar idea was used  to show the existence of auto-parallel affine submanifolds in \cite{PR2} and the existence of spherically bent submanifolds in \cite{PR1}.

In \cite{Yu}, by observing the fact that local isometries will also preserve developments of curves. We obtained an alternative form and proof of the Cartan-Ambrose-Hicks theorem by using developments of curves. In this paper, we will adapt the argument to the positive codimensional case of Cartan-Ambrose-Hicks theorem. The main difficulty is to suitably extend the notion of developments of curves to the  positive codimensional case. In this paper, we overcome this difficulty by introducing the notion of developments of curves with respect to a symmetric tensor.

Let's recall the notion of developments of curves first.
\begin{defn}\label{def-dev}
Let $(M,g)$ be a Riemannian manifold, $p\in M$ and $v:[0,T]\to T_pM$ be a smooth curve in $T_pM$. The development of $v$ is a curve $\gamma:[0,T]\to M$ such that
$$\left\{\begin{array}{ll}\gamma(0)=p&\\
\gamma'(t)=P_0^t(\gamma)(v(t))& t\in [0,T].
\end{array}\right.$$
Here $P_{t_1}^{t_2}(\gamma)$ is the parallel displacement from $\gamma(t_1)$ to $\gamma(t_2)$ along $\gamma$.
\end{defn}

In the monograph \cite{KN} by Kobayashi and Normizu, this notion is defined reversely. That is, the curve $v(t)$ in Definition \ref{def-dev} is called the development of $\gamma(t)$ in \cite{KN}. A proof of the uniqueness and local existence for developments of curves can be found in \cite{Yu}. The global existence of developments of curves is equivalent to the completeness of the Riemannian metric (see \cite{KN}). We will denote the development of $v$ as $\dev(p,v)$ when it exists. By definition \ref{def-dev}, if for any smooth curve $\gamma:[0,1]\to M$ with $\gamma(0)=p$, define
$$v_\gamma(t):=P_t^0(\gamma)(\gamma'(t)),$$
then $\gamma=\dev(p,v_\gamma).$

The alternative form of the Cartan-Ambrose-Hicks theorem we obtained in \cite{Yu} is as follows.
\begin{thm}[An alternative form of Cartan-Ambrose-Hicks theorem in \cite{Yu}]\label{thm-CAH}
Let $(M^n,g)$ and $(\wt M^n,\wt g)$ be two Riemannian manifolds with $M$ simply connected and $\wt M$ complete. Let $p\in M$, $\wt p\in \wt M$ and $\varphi:T_pM\to T_{\wt p}\wt M$ be a linear isometry.  For any smooth curve $\gamma:[0,1]\to M$ with $\gamma(0)=p$, let $\wt v_{\gamma}=\varphi(v_\gamma)$ and $\wt\gamma=\dev(\wt p,\wt v_\gamma)$.
Suppose that $R=\tau^{*}_\gamma \wt R$
for any smooth  curve $\gamma:[0,1]\to M$ where
$$\tau_\gamma=P_0^1(\wt \gamma)\circ\varphi\circ P_1^0(\gamma):T_{\gamma(1)}M\to T_{\wt\gamma(1)}\wt M.$$
Then, the map $f(\gamma(1))=\wt\gamma(1)$ from $M$ to $\wt M$ is well defined and $f$ is the local isometry from $M$ to $\wt M$ with $f(p)=\wt p$ and $f_{*p}=\varphi$.
\end{thm}
Motivated by Theorem \ref{thm-CAH}, one can use $\tau_\gamma^*\wt R$ to deal with the problem that $\wt R$ is not intrinsic when considering isometric immersions into general ambient spaces. This is a key observation of this work.

The final key step and most difficult part of this work is to extend
the notion of developments of curves to the positive codimensional case. Note that the original definition for developments of curves does not work for the positive codimensional case because developments of curves in a proper submanifold is different with that in the ambient space unless the submanifold is totally geodesic. Intuitively, for developments of curves in the positive codimensional case, we want to recover parallel displacements on submanifolds and on the normal vector bundles intrinsically by just using the second fundamental form. Due to this consideration, we introduce the following notion of developments of curves with respect to symmetric tensors. 
\begin{defn}[Development of curve w.r.t. a symmetric tensor]\label{def-g-dev}
Let $(\wt M^{n+s},\wt g)$ be a Riemannian manifold and $\wt p\in \wt M$. Let $T_{\wt p}\wt M=T^n\oplus N^s$ be an orthogonal decomposition,  $\wt h(t):[0,l]\to \Hom(T\odot T, N)$ and $\wt v:[0,l]\to T$ be smooth maps. A curve $\wt \gamma:[0,l]\to \wt M$ is called a development of $\wt v$ w.r.t. $\wt h$ if there exists a moving orthonormal frame $\left\{\wt E_A|A=1,2,\cdots,n+s\right\}$ along $\wt\gamma $ satisfying the following equations:
\begin{equation}\label{eq-g-dev}
\left\{\begin{array}{ll}
\wt\gamma(0)=\wt p\\
\wt e_a:=\wt E_a(0)\in T&a=1,2,\cdots,n\\
\wt e_\alpha:=\wt E_\alpha(0)\in N&\alpha=n+1,n+2,\cdots,n+s\\
\wt \nabla_{\wt\gamma'(t)}\wt E_a=\displaystyle\sum_{\alpha=n+1}^{n+s}\vv<\wt h(t)(\wt v(t),\wt e_a),\wt e_\alpha>\wt E_\alpha& a=1,2,\cdots,n\\
\wt\nabla_{\wt\gamma'(t)}\wt E_\alpha=-\displaystyle \sum_{a=1}^n\vv<\wt h(t)(\wt v(t),\wt e_a),\wt e_\alpha>\wt E_a&\alpha=n+1,\cdots,n+s\\
\wt\gamma'(t)=\displaystyle\sum_{a=1}^n\vv<\wt v(t),\wt e_a>\wt E_a.\\
\end{array}\right.
\end{equation}
Here $T\odot T$ means the symmetric product of $T$. Moreover, define the map $D_{t_1}^{t_2}(\wt\gamma)$ as
 $$D_{t_1}^{t_2}(\wt \gamma):T_{\wt \gamma(t_1)}\wt M\to T_{\wt \gamma(t_2)}\wt M,\ 
 \sum_{A=1}^{n+s}c_AE_A(t_1)\mapsto\sum_{A=1}^{n+s}c_A\wt E_A(t_2)$$
 for any $c_1,c_2\cdots,c_{n+s}\in \R$.
\end{defn}
It is not hard to see that the definition above is independent of the choices of the orthonormal basis $\wt e_1,\wt e_2,\cdots,\wt e_{n+s}$. Moreover, $D_{t_1}^{t_2}(\wt\gamma)$ preserves metrics by definition. When $\wt h=0$, one can see that the developments of curves w.r.t. $h$ are just the same as the classical developments of curves in Definition \ref{def-dev}, and in this case, $D_0^t(\wt \gamma)=P_0^t(\wt\gamma)$. By a similar argument as in \cite{Yu}, one can derive the equations for developments of curves w.r.t. symmetric tensors and  show its local existence, uniqueness,  and when $\wt M$ is complete, its global existence (see \ref{thm-g-dev}). We will then denote the curve $\wt \gamma$ in  Definition \ref{def-g-dev} as $\dev(\wt p,\wt v,\wt h)$ when it exists.

Before stating the main results of this paper, we first fix some notations that will be used. Let $(M^n,g)$ be a Riemannian manifold,  $(V^s,\mathfrak{h},D)$ be a smooth vector bundle on $M$ equipped with the Riemannian metric $\mathfrak{h}$ and the compatible connection $D$. Let $h\in \Gamma(\Hom(TM\odot TM,V))$.  For any $\xi\in V_p$, we will define $A_\xi:T_pM\to T_pM$ by
\begin{equation*}
\vv<A_\xi(X),Y>_g=\vv<h(X,Y),\xi>_{\mathfrak{h}}, \forall\ X,Y\in T_p(M).
\end{equation*}
We are now ready to state the first main result of this paper, the existence of isometric immersions into general ambient spaces with prescribed second fundamental form.
\begin{thm}\label{thm-main}
Let $(M^n,g)$ and $(\wt M^{n+s},\wt g)$ be two Riemannian manifolds with $M$ simply connected and $\wt M$ complete. Let $(V^s,\mathfrak{h}, D)$ be a smooth vector bundle on $M$ equipped with a Riemannian metric $\mathfrak{h}$ and a compatible connection $D$, and $h\in \Gamma(\Hom(TM\odot TM, V))$. Let $p\in M$, $\wt p\in\wt M$,  $\varphi:T_pM\oplus V_p\to T_{\wt p}\wt M$ be a linear isometry, and let $T=\varphi(T_pM)$ and $N=\varphi(V_p)$. For any smooth curve $\gamma:[0,1]\to M$ with $\gamma(0)=p$, let $\wt v_\gamma(t)=\varphi(v_\gamma(t))$ and $\wt \gamma=\dev(\wt p, \wt v_\gamma, \wt h_\gamma)$ where
$$\wt h_\gamma(t)=(\varphi^{-1})^*P_t^0(\gamma)h,\ \forall t\in [0,1].$$
Moreover, for any smooth curve $\gamma:[0,1]\to M$ with $\gamma(0)=p$, suppose that\\
(1) for any $X,Y,Z,W\in T_{\gamma(1)}M$,
$$R(X,Y,Z,W)=(\tau_\gamma^*\wt R)(X,Y,Z,W)+\vv<h(X,W),h(Y,Z)>-\vv<h(X,Z),h(Y,W)>,$$
where $R$ and $\wt R$ are curvature tensors of $M$ and $\wt M$ respectively;\\
(2)  for any tangent vectors $X,Y,Z\in T_{\gamma(1)}M$ and $\xi\in V_{\gamma(1)}$,
$$\vv<(D_{X}h)(Y,Z)-(D_Yh)(X,Z),\xi>=(\tau^*_\gamma\wt R)(Z,\xi,X,Y);$$
(3) for any $X,Y\in T_{\gamma(1)}M$ and $\xi,\eta\in V_{\gamma(1)}$,
$$R^V(\xi,\eta,X,Y)=(\tau_{\gamma}^*\wt R)(\xi,\eta,X,Y)+\vv<A_\xi(Y),A_\eta(X)>-\vv<A_\eta(Y),A_\xi(X)>,$$
where $R^V$ is the curvature tensor of the vector bundle $V$. Here
    $$\tau_\gamma=D_0^1(\wt \gamma)\circ\varphi\circ P_1^0(\gamma):T_{\gamma(1)}M\oplus V_{\gamma(1)}\to T_{\wt\gamma(1)}\wt M.$$
Then, we have the following conclusions.
\begin{enumerate}
\item The map  $f(\gamma(1))=\wt\gamma(1)$ from $M$ to $\wt M$ is well defined and is an isometric immersion from $M$ to $\wt M$.
\item The map $\wt f:V\to T^\perp M$ with $\wt f|_{\gamma(1)}=\tau_\gamma|_{V_{\gamma(1)}}$ is well defined, where $T^\perp M$ is the normal bundle of the immersion $f:M\to \wt M$. Moreover, $\wt f$ preserves metrics and connections.
\item $\wt f^*h_{\wt M}=h$ where $h_{\wt M}$ is the second fundamental form of the isometric immersion $f:M\to \wt M$.
\end{enumerate}
\end{thm}
When the ambient space $(\wt M,\wt g)$ in Theorem \ref{thm-main} is a space form with constant sectional curvature $K$,
$$\tau_\gamma^*\wt R(X,Y,Z,W)=K\left(\vv<X,W>\vv<Y,Z>-\vv<X,Z>\vv<Y,W>\right)$$
because $\tau_\gamma$ is a linear isometry. It is then clear that Theorem \ref{thm-main} contains the existence part of the fundamental theorem for submanifolds in space forms (\cite{Bo,Sz,Te}) because the requirements (1), (2) and (3) in Theorem \ref{thm-main} are independent of the curve $\gamma$ in this case.  One advantage of Theorem \ref{thm-main} is that it also provides a geometric construction of the isometric immersion when it exists. Moreover, when $s=0$, Theorem \ref{thm-main} reduces to Theorem \ref{thm-CAH}. So, Theorem \ref{thm-main} is in fact a Cartan-Ambrose-Hicks theorem for isometric immersions.

Finally, we would like to mention that the local version of the Cartan-Ambrose-Hicks theorem is the Cartan isometry theorem in Riemannian geometry (see \cite{CE}) which only requires weaker assumptions than the global version. So, we also provide the local version of Theorem \ref{thm-main} below.
\begin{thm}[Local version of Theorem \ref{thm-main}]\label{thm-Cartan}
Let $(M^n,g)$ be a Riemannian manifold and $p\in M$. Let $\widehat\Omega\subset T_pM$ be a star-shaped domain centered at $0$ such that $\exp_p:\widehat\Omega\to \Omega:=\exp_p(\widehat\Omega)$ is a diffeomorphism.  Let $(\wt M^{n+s},\wt g)$ be a complete Riemannian manifold and $\wt p\in \wt M$. Let $(V^s,\mathfrak{h}, D)$ be a smooth vector bundle on $\Omega$ equipped with a Riemannian metric $\mathfrak{h}$ and a compatible connection $D$, and $h\in \Gamma(\Hom(T\Omega\odot T\Omega, V))$. Let $p\in M$, $\wt p\in\wt M$,  $\varphi:T_pM\oplus V_p\to T_{\wt p}\wt M$ be a linear isometry, and let $T:=\varphi(T_pM)$ and $N:=\varphi(V_p)$. For any geodesic $\gamma:[0,1]\to \Omega$ with $\gamma(0)=p$, let $\wt v_\gamma(t)=\varphi(\gamma'(0))$ and $\wt \gamma=\dev(\wt p, \wt v_\gamma, \wt h_\gamma)$ where
$$\wt h_\gamma(t)=(\varphi^{-1})^*P_t^0(\gamma)h$$
for $t\in [0,1]$. Moreover, for any geodesic $\gamma:[0,1]\to \Omega$ with $\gamma(0)=p$, suppose that\\
(1) for any $X,Y,Z,W\in T_{\gamma(1)}M$,
$$R(X,Y,Z,W)=(\tau_\gamma^*\wt R)(X,Y,Z,W)+\vv<h(X,W),h(Y,Z)>-\vv<h(X,Z),h(Y,W)>,$$
where $R$ and $\wt R$ are curvature tensors of $M$ and $\wt M$ respectively;\\
(2)  for any tangent vectors $X,Y,Z\in T_{\gamma(1)}M$ and $\xi\in V_{\gamma(1)}$,
$$\vv<(D_{X}h)(Y,Z)-(D_Yh)(X,Z),\xi>=(\tau^*_\gamma\wt R)(Z,\xi,X,Y);$$
(3) for any $X,Y\in T_{\gamma(1)}M$ and $\xi,\eta\in V_{\gamma(1)}$,
$$R^V(\xi,\eta,X,Y)=(\tau_{\gamma}^*\wt R)(\xi,\eta,X,Y)+\vv<A_\xi(Y),A_\eta(X)>-\vv<A_\eta(Y),A_\xi(X)>,$$
where $R^V$ is the curvature tensor of the vector bundle $V$. Here
    $$\tau_\gamma=D_0^1(\wt \gamma)\circ\varphi\circ P_1^0(\gamma):T_{\gamma(1)}M\oplus V_{\gamma(1)}\to T_{\wt\gamma(1)}\wt M.$$
Then, we have the following conclusions.
\begin{enumerate}
\item The map  $f(\gamma(1))=\wt\gamma(1)$ from $\Omega$ to $\wt M$ is an isometric immersion from $\Omega$ to $\wt M$.
\item The bundle map $\wt f:V\to T^\perp \Omega$ with $\wt f|_{\gamma(1)}=\tau_\gamma|_{V_{\gamma(1)}}$ preserves metrics and connections where $T^\perp \Omega$ is the normal bundle of the immersion $f:\Omega\to \wt M$.
\item $\wt f^*h_{\wt M}=h$ where $h_{\wt M}$ is the second fundamental form of the isometric immersion $f:\Omega\to \wt M$.
\end{enumerate}
\end{thm} 
The main difference between the local version (Theorem \ref{thm-Cartan}) and the global version (Theorem \ref{thm-main}) is that one only need to verify the Gauss-Codazzi-Ricci equations for geodesics in $\Omega$ starting at $p$ which is uniquely determined by their terminal points.

In this paper, for simplicity, we adopt the following conventions of notations:
\begin{enumerate}
\item Capital Latin letters such as $A,B,C,D$ denote indices in $\{1,2,\cdots,n+s\}$;
\item Latin letters in lower case such as $a,b,c,d$ denote indices in $\{1,2,\cdots,n\}$;
\item Greek letters such as $\alpha,\beta$ denote indices in $\{n+1,n+2,\cdots,n+s\}$;
\item Repeated indices indicate summation;  
\item $u'$ means taking derivative of $u$ w.r.t. $t$. 
\end{enumerate}

The rest of the paper is organized as follows. In Section 2, we show the local existence, uniqueness and global existence of developments of curves w.r.t. a symmetric tensor and derive the equation for the variation field for a family of developments of curves w.r.t. symmetric tensors. In Section 3, we prove Theorem \ref{thm-Cartan} and Theorem \ref{thm-main}. 
\section{Developments of curves with respect to symmetric tensors}
In this section, we will show the local existence, uniqueness and global existence of the developments of curves w.r.t. symmetric tensors, derive the equations for the variation field for a family of  developments of curves w.r.t.  symmetric tensors.
\begin{thm}\label{thm-g-dev} Let the notations be the same as in Definition \ref{def-g-dev}. Then, $\wt \gamma$  is unique and exists for a short time. Moreover, if the Riemannian manifold $(\wt M,\wt g)$ is complete, then $\wt \gamma$ exists all over $[0,l]$.
\end{thm}
\begin{proof}
Let $(x_1,x_2,\cdots,x_{n+s})$ be a local coordinate at $\wt p$ with $x_A(\wt p)=0$ and $\frac{\p}{\p x_A}(p)=\wt e_A$ for $A=1,2,\cdots, n+s$, where $\wt e_1,\cdots,\wt e_n$ is an orthonormal basis of $T$ and $\wt e_{n+1},\cdots,\wt e_{n+s}$ is an orthonormal basis of $N$. 

Suppose that $\wt \gamma(t)=(x_1(t),x_2(t),\cdots,x_{n+s}(t))$ in the local coordinate, $\wt E_{A}=X_{AB}\frac{\p}{\p x_B},$ $\wt v=v_a\wt e_a$ and
\begin{equation*}
h_{ab}^\alpha(t)=\vv<\wt h(t)(\wt e_a,\wt e_b),\wt e_\alpha>.
\end{equation*}
Substituting all the above into \eqref{eq-g-dev}, we have
\begin{equation}
\left\{\begin{array}{l} X'_{aA}+X_{aB}X_{bC}v_b\wt \Gamma_{BC}^A-h_{ab}^\alpha v_bX_{\alpha A}=0\\
X'_{\alpha A}+X_{\alpha B}X_{bC}v_b\wt \Gamma_{BC}^A+ h_{ab}^\alpha v_b X_{a A}=0\\
x'_A-v_aX_{aA}=0\\
x_A(0)=0\\
X_{AB}(0)=\delta_{AB}.
\end{array}\right.
\end{equation}
Here $\wt \Gamma_{AB}^C$'s are the Christoffel symbols for $\wt g$. By standard theory for ODEs, we get the local existence and uniqueness for the solution of the equation. 

Next we show that $\wt E_1,\cdots,\wt E_{n+s}$ is an orthonormal frame along $\wt \gamma$. In fact, let $Y_{AB}=\vv<\wt E_A,\wt E_B>-\delta_{AB}$. Then, $Y_{AB}(0)=0$ for any $A,B=1,2,\cdots,n+s$. Moreover, by \eqref{eq-g-dev},
\begin{equation*}
\begin{split}
Y'_{a\alpha}=&\frac{d}{dt}\vv<\wt E_a,\wt E_\alpha>\\
=&\vv<\wt h(t)(\tilde v(t),\wt e_a),\wt e_\beta>\vv<\wt E_\beta,\wt E_\alpha>-\vv<\wt h(t)(\wt v(t),\wt e_b),\wt e_\alpha>\vv<\wt E_b,\wt E_a>\\
=&\vv<\wt h(t)(\wt v(t),\wt e_a),\wt e_\beta>Y_{\beta\alpha}-\vv<\wt h(t)(\wt v(t),\wt e_b),\wt e_\alpha>Y_{ba},
\end{split}
\end{equation*}
and similarly,
\begin{equation*}
Y'_{ab}=\vv<\wt h(t)(\wt v(t),\wt e_a),\wt e_\alpha>Y_{\alpha b}+\vv<\wt h(t)(\wt v(t),\wt e_b),\wt e_\alpha>Y_{\alpha a}
\end{equation*}
and
\begin{equation*}
Y'_{\alpha\beta}=-\vv<\wt h(t)(\wt v(t),\wt e_a),\wt e_\alpha>Y_{a\beta}-\vv<\wt h(t)(\wt v(t),\wt e_a),\wt e_\beta>Y_{a\alpha}.
\end{equation*}
So, $Y$ satisfies a first order homogeneous linear system of ODEs with initial data $Y_{AB}(0)=0$. This implies that $Y_{AB}(t)=0$ for any $t$ and thus 
$\wt E_1,\cdots,\wt E_{n+s}$ is an orthonormal frame along $\wt\gamma$. This completes the proof of local existence and uniqueness for developments of curves w.r.t. symmetric tensors.

Finally, we show the global existence of $\wt\gamma$ when $(\wt M,\wt g)$ is complete. In fact, if $\wt \gamma$ does not exist on $[0,l]$, let $[0,T)$ with $T\leq l$ be the maximum existence interval for $\wt\gamma$. Then, for any $t_1,t_2\in [0,T)$ with $t_1<t_2$, we have
\begin{equation*}
d(\wt\gamma(t_1),\wt\gamma(t_2))\leq \int_{t_1}^{t_2}\|\wt\gamma'\|dt=\int_{t_1}^{t_2}\|v(t)\|dt\leq L|t_2-t_1|
\end{equation*}
where $L=\max_{t\in [0,l]}\|v(t)\|$. So, by completeness of $\wt M$, we know that $\lim_{t\to T^{-}}\wt\gamma(t)$ exists and $\wt\gamma$ can be further extended to $[0,T]$. This is a contradiction.
\end{proof}
Next, we come to derive the equations for the variation field of a family of developments of curves w.r.t. symmetric tensors.
\begin{thm}\label{thm-g-Jacobi}
Let $(\wt M^{n+s},\wt g)$ be a Riemannian manifold and $\wt p\in \wt M$. Let $T_{\wt p}\wt M=T^n\oplus N^s$ be an orthogonal decomposition, and $\wt v(u,t):I\times [0,1]\to T$ and $\wt h(u,t):I\times[0,1]\to \Hom(T\odot T, N)$ be smooth maps where $I$ is some interval. Let
\begin{equation*}
\wt \Phi(u,t)=\wt \gamma_u(t):=\dev\left(\wt p, \wt v(u,\cdot),\wt h(u,\cdot)\right)(t),
\end{equation*}
$\wt e_1,\cdots,\wt e_n$ be an orthonormal basis for $T$, and $\wt e_{n+1},\cdots,\wt e_{n+s}$ be an orthonormal basis for $N$. Moreover, let $\wt E_A(u,t)=D_0^t(\wt\gamma_u)(\wt e_A)$,
$\wt v(u,t)=v_a(u,t)\wt e_a$ and
\begin{equation*}
\wt h(u,t)(\wt e_a,\wt e_b)=h_{ab}^\alpha(u,t)\wt e_\alpha.
\end{equation*}
Suppose that 
\begin{equation*}
\frac{\p\wt \Phi}{\p u}=\wt U_A\wt E_A\ \mbox{and }\wt \nabla_{\pu}\wt E_A=\wt X_{AB}\wt E_B.
\end{equation*}
Then, $\wt X_{AB}=-\wt X_{BA}$, and
\begin{equation}\label{eq-g-Jacobi-2}
\left\{\begin{array}{rl}\wt U''_a=&2\wt U'_\alpha h_{ab}^\alpha v_b+\wt U_\alpha\p_t(h_{ab}^\alpha v_b)+\wt U_ch_{cd}^\alpha h_{ab}^\alpha v_d v_b+\wt R_{bacA} \wt U_A v_bv_c\\
&+(\p_t v_b)\wt X_{ba}-v_b v_c h_{bc}^\alpha\wt X_{a\alpha }+\p_u\p_tv_a\\
\wt U''_\alpha=&-2\wt U'_a h_{ab}^\alpha v_b-\wt U_{c}\p_t( h_{bc}^\alpha v_b)+
\wt U_\beta h_{bc}^\beta h_{ab}^\alpha  v_a v_c+\wt R_{b\alpha aA}\wt U_A v_av_b\\
&+(\p_t v_a)\wt X_{a\alpha}+v_av_bh_{ab}^\beta\wt X_{\beta\alpha}+\p_u( v_av_bh_{ab}^\alpha)\\
\wt X'_{ab}=&\wt X_{a\alpha}h_{bc}^\alpha v_c-h_{ac}^\alpha v_c\wt X_{b\alpha}+\wt R_{abcA}\wt U_A v_c\\
\wt X'_{a\alpha}=&-\wt X_{ab}h_{bc}^\alpha v_c+ h_{ab}^\beta v_b\wt X_{\beta\alpha}+\wt R_{a\alpha bA}\wt U_Av_b+\p_u(h_{ab}^\alpha v_b)\\
\wt X'_{\alpha\beta}=&\wt X_{a\alpha} h_{ab}^\beta v_b-\wt X_{a\beta}h_{ab}^\alpha v_b+\wt R_{\alpha\beta aA}v_a\wt U_A\\
\wt U_A(u,0)=&\wt X_{AB}(u,0)=\wt U'_{\alpha}(u,0)=0\\
\wt U'_a(u,0)=&\p_uv_a(u,0).\\
\end{array}\right.
\end{equation}
Here $\wt R_{ABCD}=\wt R(\wt E_A,\wt E_B,\wt E_C,\wt E_D)$ with $\wt R$ the curvature tensor of $(\wt M,\wt g)$.
\end{thm}
\begin{proof} By that $\vv<\wt E_A,\wt E_B>=\delta_{AB}$, we have
\begin{equation*}
\wt X_{AB}=\vv<\wt\nabla_\pu \wt E_A,\wt E_B>=-\vv<\wt E_A,\wt\nabla_\pu \wt E_B>=-\wt X_{BA}.
\end{equation*}
By \eqref{eq-g-dev}, we have
\begin{equation}\label{eq-t}
\frac{\p\wt \Phi}{\p t}=\wt\gamma'_u=v_a(u,t)\wt E_a(u,t).
\end{equation}
So, by \eqref{eq-g-dev} again,
\begin{equation}\label{eq-1-t}
\begin{split}
\wt \nabla_{\pt}\pt=\p_t v_a\wt E_a+v_a\wt\nabla_{\pt}\wt E_a=\p_t v_a\wt E_a+v_a v_bh_{ab}^\alpha\wt E_\alpha,\\
\end{split}
\end{equation}
and
\begin{equation}\label{eq-1-u}
\begin{split}
\wt \nabla_{\pt}\pu=&\wt U'_A\wt E_A+\wt U_a\wt\nabla_{\pt}\wt E_a+\wt U_\alpha \wt \nabla_{\pt}\wt E_\alpha\\
=&\wt U'_A\wt E_A+\wt U_a h_{ab}^\alpha v_b\wt E_\alpha-\wt U_\alpha  h_{ab}^\alpha v_b\wt E_a.\\
\end{split}
\end{equation}
Then,
\begin{equation}\label{eq-2-u-1}
\begin{split}
&\wt\nabla_{\pt}\wt\nabla_{\pt}\pu\\
=&\wt\nabla_{\pt}\left(\wt U'_A\wt E_A+\wt U_a h_{ab}^\alpha v_b\wt E_\alpha-\wt U_\alpha  h_{ab}^\alpha v_b\wt E_a\right)\\
=&\wt U''_A\wt E_A+2\wt U'_a h_{ab}^\alpha v_b\wt E_\alpha-2\wt U'_\alpha  h_{ab}^\alpha v_b\wt E_a\\
&+\left(\wt U_{c}\p_t( h_{bc}^\alpha v_b)-
\wt U_\beta h_{bc}^\beta h_{ab}^\alpha  v_av_c\right)\wt E_\alpha-\left(\wt U_\alpha\p_t\left( h_{ab}^\alpha v_b\right)+\wt U_c h_{cd}^\alpha h_{ab}^\alpha v_d v_b\right)\wt E_a.
\end{split}
\end{equation}
On the other hand, by \eqref{eq-1-t},
\begin{equation}\label{eq-2-u-2}
\begin{split}
&\wt\nabla_{\pt}\wt \nabla_{\pt}\pu\\
=&\wt\nabla_{\pt}\wt\nabla_{\pu}\pt\\
=&\wt\nabla_{\pu}\wt\nabla_{\pt}\pt+\wt R\left(\pt,\pu\right)\pt\\
=&\wt\nabla_{\pu}\left(\p_tv_a\wt E_a+v_a v_bh_{ab}^\alpha\wt E_\alpha\right)+\wt R_{bBaA}\wt U_Av_av_b\wt E_B\\
=&\left(\wt R_{bacA}\wt U_Av_bv_c+\p_u\p_tv_a+(\p_t v_b)\wt X_{ba}+ v_b v_ch_{bc}^\alpha\wt X_{\alpha a}\right)\wt E_a+\\
&\left(\wt R_{b\alpha aA}\wt U_Av_a v_b+(\p_tv_a)\wt X_{a\alpha}+\p_u( v_av_bh_{ab}^\alpha)+v_av_bh_{ab}^\beta\wt X_{\beta\alpha}\right)\wt E_\alpha.
\end{split}
\end{equation}
Comparing \eqref{eq-2-u-1} and \eqref{eq-2-u-2}, and noting that $\wt X_{AB}=-\wt X_{BA}$, we get the first two equations of \eqref{eq-g-Jacobi-2}.

Moreover, note that
\begin{equation}\label{eq-dt-i-1}
\begin{split}
\wt\nabla_{\pt}\wt\nabla_{\pu}\wt E_a=&\wt \nabla_{\pt}(\wt X_{a\alpha}\wt E_\alpha)+\wt \nabla_{\pt}(\wt X_{ab}\wt E_b)\\
=&\left(\wt X'_{a\alpha}+\wt X_{ab}h_{bc}^\alpha v_c\right)\wt E_\alpha+\left(\wt X'_{ab}-\wt X_{a\alpha} h_{bc}^\alpha v_c\right)\wt E_b
\end{split}
\end{equation}
and on the other hand, by \eqref{eq-g-dev},
\begin{equation}\label{eq-dt-i-2}
\begin{split}
&\wt\nabla_{\pt}\wt\nabla_{\pu}\wt E_a\\
=&\wt\nabla_{\pu}\wt\nabla_{\pt}\wt E_a+\wt R\left(\pt,\pu\right)\wt E_a\\
=&\left( h_{ac}^\alpha v_c\wt X_{\alpha b}+\wt R_{abcA}\wt U_A v_c\right)\wt E_b+\left(\p_u(h_{ab}^\alpha v_b)+h_{ab}^\beta v_b\wt X_{\beta\alpha}+\wt R_{a\alpha bA}\wt U_A v_b\right)\wt E_\alpha.
\end{split}
\end{equation}
Hence, by comparing \eqref{eq-dt-i-1} and \eqref{eq-dt-i-2}, and noting that $\wt X_{AB}=-\wt X_{BA}$, we get the third and the fourth equations of \eqref{eq-g-Jacobi-2}. Similarly, by computing $\wt\nabla_\pt\wt\nabla_\pu\wt E_\alpha$ in the two ways as above, one get the fifth equation of \eqref{eq-g-Jacobi-2}.

Finally, $\wt U_A(u,0)=0$ because $\Phi(u,0)=\wt p$, and
$X_{AB}(u,0)=0$ because $\wt E_A(u,0)=\wt e_A$ is independent of $u$. Moreover, by \eqref{eq-t}, we have
\begin{equation*}
\begin{split}
\wt U'_A(u,0)\wt e_A=\wt\nabla_{\pt}\pu\bigg|_{t=0}
=\wt\nabla_{\pu}\pt\bigg|_{t=0}
=\wt\nabla_{\pu}(v_a\wt E_a)\bigg|_{t=0}
=(\p_uv_a)(u,0)\wt e_a.
\end{split}
\end{equation*}
So, $\wt U'_a(u,0)=\p_uv_a(u,0)$ and $\wt U'_\alpha(u,0)=0$ and we have obtained all the initial data in \eqref{eq-g-Jacobi-2}. This completes the proof of the theorem.
\end{proof}
\section{Existence of isometric immersions with prescribed second fundamental form}
In this section, by using a similar idea to that in  \cite{Yu}, we prove Theorem \ref{thm-main} and Theorem \ref{thm-Cartan}, the existence of isometric immersions into general ambient spaces assuming the Gauss-Codazzi-Ricci equations are satisfied.  Because the proof of the  local existence of isometric immersion is simpler, we first prove Theorem \ref{thm-Cartan}.

\begin{proof}[Proof of Theorem \ref{thm-Cartan}]
Let $e_1,e_2,\cdots,e_n$ and $e_{n+1},e_{n+2},\cdots,e_{n+s}$ be orthonormal basis of $T_pM$ and $V_p$ respectively, and let $\wt e_A=\varphi(e_A)$ for $A=1,2,\cdots,n+r$. Let $\Phi(u,t):(-\e,\e)\times [0,1]\to \Omega$ be a smooth family of geodesics starting at $p$. That is, for each $u\in (-\e,\e)$, $\gamma_u(\cdot):=\Phi(u,\cdot)$ is a geodesic in $\Omega$ with $\gamma_u(0)=p$.
Let $\wt\Phi(u,t)=\wt\gamma_u(t)$. Moreover, let $E_A(u,t)=P_0^t(\gamma_u)(e_A)$ and 
$\wt E_A(u,t)=D_0^t(\wt\gamma_u)(\wt e_A)$ for $A=1,2,\cdots, n+s$.

Let $\frac{\p\Phi}{\p u}=U_a E_a$, $\nabla_{\frac{\p}{\p u}}E_{a}=X_{ab}E_b$ and 
$D_{\frac{\p}{\p u}}E_\alpha=X_{\alpha\beta}E_\beta$. Suppose that 
$$\gamma_u'(0)=v_a(u)e_a.$$ 
Then, by the same computation as in the proof of Theorem \ref{thm-g-Jacobi}, we know that $U$ and  $X$ satisfy the following system of ODEs:
\begin{equation}\label{eq-Jacobi-field}
\begin{split}
\left\{\begin{array}{l}U''_a=R_{cadb}v_cv_dU_b\\
X'_{ab}=R_{abdc}v_dU_c\\
X'_{\alpha\beta}=R^V_{\alpha\beta ab}v_aU_b\\
U_a(u,0)=X_{ab}(u,0)=X_{\alpha\beta}(u,0)=0\\
U'_a(u,0)=\p_uv_a(u),
\end{array}\right.
\end{split}
\end{equation}
where $R_{abcd}=R(E_a,E_b,E_c,E_d)$ and $R^V_{\alpha\beta ab}=R^V(E_\alpha,E_\beta,E_a,E_b)$.

On the other hand, let
\begin{equation*}
\frac{\p\wt\Phi}{\p u}=\wt U_AE_A\mbox{ and } \wt \nabla_{\pu}\wt E_A=\wt X_{AB}\wt E_B.
\end{equation*}
By Theorem \ref{thm-g-Jacobi}, $\wt U$ and $\wt X$ satisfy \eqref{eq-g-Jacobi-2}. We claim that
\begin{equation}\label{eq-solution-1}
\left\{\begin{array}{l}\wt U_a=U_a\\
\wt U_\alpha=0\\
\wt X_{ab}=X_{ab}\\
\wt X_{\alpha\beta}=X_{\alpha\beta}\\
\wt X_{a\alpha}=h_{ab}^\alpha U_b.
\end{array}\right.
\end{equation}
Here 
$$h_{ab}^\alpha=\vv<h(E_a,E_b),E_\alpha>=\vv<\wt h(\wt e_a,\wt e_b),\wt e_\alpha>.$$
 By Theorem \ref{thm-g-Jacobi}, we only need to verify that the $\wt U_A$'s and $\wt X_{AB}$'s given above in \eqref{eq-solution-1} satisfy the Cauchy problem \eqref{eq-g-Jacobi-2}. 
 
 First of all, by the initial data in \eqref{eq-Jacobi-field}, the initial data of $\wt U$ and $\wt X$ given by \eqref{eq-solution-1} satisfy the initial data in \eqref{eq-g-Jacobi-2}. We now verify that $\wt U$ and $\wt X$ given by \eqref{eq-solution-1} satisfy the equations in \eqref{eq-g-Jacobi-2}.
 
By \eqref{eq-Jacobi-field}, \eqref{eq-solution-1} and assumption (1) of Theorem \ref{thm-main},
\begin{equation*}
\begin{split}
\wt U''_a=&U''_a\\
=&R_{cadb}v_cv_dU_b\\
=&\left(\wt R_{cadb}+h^\alpha_{cb}h_{ad}^\alpha-h^\alpha_{cd}h_{ab}^\alpha\right)v_cv_dU_b\\
=&2\wt U'_\alpha h_{ab}^\alpha v_b+\wt U_\alpha\p_t(h_{ab}^\alpha v_b)+\wt U_ch_{cd}^\alpha h_{ab}^\alpha v_d v_b+\wt R_{bacA} \tilde U_A v_bv_c\\
&+(\p_t v_b)\wt X_{ba}-v_b v_c h_{bc}^\alpha\wt X_{a\alpha }+\p_u\p_tv_a\\
\end{split}
\end{equation*}
by noting that $\wt U_\alpha=0$ and $\p_tv_a=0$. So, the first equation in \eqref{eq-g-Jacobi-2} is satisfied by \eqref{eq-solution-1}.

By assumption (2) of Theorem \ref{thm-main},
\begin{equation*}
\begin{split}
&\p_u(v_av_b h_{ab}^\alpha)\\
=&\p_u\vv<h\left(\pt,\pt\right),E_\alpha>\\
=&\vv<D_\pu \left(h\left(\pt,\pt\right)\right),E_\alpha>+\vv<h\left(\pt,\pt\right),D_\pu E_\alpha>\\
=&\vv<\left(D_\pu h\right)\left(\pt,\pt\right),E_\alpha>+2\vv<h\left(\nabla_\pu \pt,\pt\right),E_\alpha>+\vv<h\left(\pt,\pt\right),D_\pu E_\alpha>\\
=&\vv<\left(\nabla_\pt h\right)\left(\pu,\pt\right),E_\alpha>+\wt R_{a\alpha cb}v_av_bU_c\\
&+2\vv<h\left(\nabla_\pt \pu,\pt\right),E_\alpha>+\vv<h\left(\pt,\pt\right),D_\pu E_\alpha>\\
=&(\p_t h_{ab}^\alpha) U_av_b+\wt R_{a\alpha cb}v_av_bU_c+2U'_ah_{ab}^\alpha v_b+h_{ab}^\beta v_av_b X_{\alpha\beta}.
\end{split}
\end{equation*}
Substituting the last equation into the RHS of the second equation in \eqref{eq-g-Jacobi-2}, we have
\begin{equation*}
\begin{split}
&-2\wt U'_a h_{ab}^\alpha v_b-\wt U_{c}\p_t( h_{bc}^\alpha v_b)+
\wt U_\beta h_{bc}^\beta h_{ab}^\alpha  v_a v_c\\
&+\wt R_{b\alpha aA}\wt U_A v_av_b+(\p_t v_a)\wt X_{a\alpha}+\p_u( v_av_bh_{ab}^\alpha)+v_av_bh_{ab}^\beta\wt X_{\beta\alpha}\\
=&-2U'_a h_{ab}^\alpha v_b-U_{c}\p_t( h_{bc}^\alpha v_b)+
 \wt R_{b\alpha ac} U_c v_av_b\\
&+(\p_t v_a)h_{ab}^\alpha U_b+\left((\p_t h_{ab}^\alpha) U_av_b+\wt R_{a\alpha cb}v_av_b U_c+2U'_ah_{ab}^\alpha v_b+h_{ab}^\beta v_av_bX_{\alpha\beta}\right)+v_av_bh_{ab}^\beta X_{\beta\alpha}\\
=&0\\
=&\wt U''_\alpha
\end{split}
\end{equation*}
by noting that $\wt U_\alpha=0$ in \eqref{eq-solution-1}, $\p_t v_a=0$ and $X_{\alpha\beta}=-X_{\beta\alpha}$. Hence, the second equation in \eqref{eq-g-Jacobi-2} is satisfied by \eqref{eq-solution-1}.

By \eqref{eq-Jacobi-field}, \eqref{eq-solution-1} and assumption (1) of Theorem \ref{thm-main},
\begin{equation*}
\begin{split}
\wt X'_{ab}=&X'_{ab}\\
=&R_{abcd}v_cU_d\\
=&\left(\wt R_{abcd}+h_{ad}^\alpha h_{bc}^\alpha-h_{ac}^\alpha h_{bd}^\alpha\right)v_cU_d\\
=&\wt X_{a\alpha}h_{bc}^\alpha v_c-h_{ac}^\alpha v_c\wt X_{b\alpha}+\wt R_{abcA}\wt U_A v_c.
\end{split}
\end{equation*}
So, the third equation of \eqref{eq-g-Jacobi-2} is satisfied by \eqref{eq-solution-1}.

By \eqref{eq-Jacobi-field}, \eqref{eq-solution-1} and assumption (2) of Theorem \ref{thm-main},
\begin{equation*}
\begin{split}
&\wt X'_{a\alpha}\\
=&\p_t(h_{ab}^\alpha U_b)\\
=&\p_t\vv<h\left(\pu,E_a\right),E_\alpha>\\
=&\vv<\left(D_\pt h\right)\left(\pu,E_a\right),E_\alpha>+\vv<h\left(\nabla_\pt\pu,E_a\right),E_\alpha>\\
=&\vv<\left(D_\pu h\right)\left(\pt,E_a\right),E_\alpha>+\wt R_{a\alpha bc} U_cv_b+\vv<h\left(\nabla_\pu\pt,E_a\right),E_\alpha>\\
=&\p_u\vv< h\left(\pt,E_a\right),E_\alpha>-\vv<h\left(\pt,\nabla_\pu E_a\right),E_\alpha>-\vv<h\left(\pt, E_a\right),\nabla_\pu E_\alpha>\\
&+\wt R_{a\alpha bA}\wt U_Av_b\\
=&\p_u(h_{ab}^\alpha v_b)-X_{ab}h_{bc}^\alpha v_c+ h_{ab}^\beta v_b X_{\beta\alpha}+\wt R_{a\alpha bA}\wt U_Av_b\\
=&-\wt X_{ab}h_{bc}^\alpha v_c+ h_{ab}^\beta v_b\wt X_{\beta\alpha}+\wt R_{a\alpha bA}\wt U_Av_b+\p_u(h_{ab}^\alpha v_b).
\end{split}
\end{equation*}
So, the fourth equation of \eqref{eq-g-Jacobi-2} is satisfied by \eqref{eq-solution-1}.

Finally, by assumption (3) of Theorem \ref{thm-main}, \eqref{eq-Jacobi-field} and \eqref{eq-solution-1},
\begin{equation*}
\begin{split}
&\wt X'_{\alpha\beta}\\
=&X'_{\alpha\beta}\\
=&R^V_{\alpha\beta ab}v_aU_b\\
=&\left(\wt R_{\alpha\beta ab}+h_{bc}^\alpha h_{ac}^\beta-h_{ac}^\alpha h_{bc}^\beta\right)v_aU_b\\
=&\wt X_{a\alpha} h_{ab}^\beta v_b-\wt X_{a\beta}h_{ab}^\alpha v_b+\wt R_{\alpha\beta aA}v_a\wt U_A.
\end{split}
\end{equation*}
So, the fifth equation of \eqref{eq-g-Jacobi-2} is satisfied by \eqref{eq-solution-1}. This completes the proof the claim.

Note that 
$$U_af_*(E_a)=f_*\left(\frac{\p\Phi}{\p u}\right)=\frac{\p\wt\Phi}{\p u}=\wt U_A\wt E_A=U_a\wt E_a.$$
So, $f_*(E_a)=\wt E_a$ and $f$ is an isometric immersion. This gives us the first conclusion of Theorem \ref{thm-Cartan}. 

Note that $\wt f(E_\alpha)=\wt E_\alpha$ by definition. So, $\wt f$ preserves metrics. Moreover, 
\begin{equation}\label{eq-conn}
\wt f\left(D_\pu E_\alpha\right)=\wt f(X_{\alpha\beta}E_\beta)=X_{\alpha\beta}\wt E_\beta.
\end{equation}
By that $f_{*}(E_a)=\wt E_a$, we have $\left(\wt X_{\alpha A}\wt E_A\right)^\perp=X_{\alpha\beta}\wt E_\beta$. So, by \eqref{eq-conn},
\begin{equation*}
\begin{split}
\wt f\left(D_\pu E_\alpha\right)=\left(\wt\nabla_{\pu}\wt E_\alpha\right)^\perp=\nabla^\perp_{\pu}\wt E_\alpha=\nabla^\perp_{f_*(\pu)}\wt f(E_\alpha).
\end{split}
\end{equation*}
This means that $\wt f$ preserves connections and the second conclusion of the theorem is proved.

Finally, note that
\begin{equation*}
\begin{split}
\wt f \left(h\left(\pu,E_b\right)\right)=\wt f(U_ah_{ab}^\alpha E_\alpha)= U_ah_{ab}^\alpha \wt E_\alpha
\end{split}
\end{equation*}
and on the other hand, by \eqref{eq-solution-1},
\begin{equation*}
\begin{split}
h_{\wt M}\left(\pu, \wt E_b\right)=\left(\wt \nabla_{\pu}\wt E_b\right)^\perp=\wt X_{b\alpha}\wt E_\alpha=U_ah_{ab}^\alpha \wt E_\alpha.
\end{split}
\end{equation*}
Thus, 
$$\wt f \left(h\left(\pu,E_b\right)\right)=h_{\wt M}\left(\pu, \wt E_b\right)$$
and hence $\wt f^*h_{\wt M}=h$. This proves the third conclusion of the theorem.
\end{proof}
\begin{rem}
\begin{enumerate}
\item Although we only assume the assumptions (1),(2),(3) in Theorem \ref{thm-Cartan} is true at $\gamma(1)$ for any geodesic $\gamma$ in $\Omega$ starting at $p$, they are also true at $\gamma(t)$ for any $t\in [0,1]$ because linearly reparameterization of a geodesic is still a geodesic and parallel displacements are independent of parametrization.
\item Although we also prescribe the normal connection in Theorem \ref{thm-Cartan}, we still call the isometric immersion an isometric immersion with prescribed second fundamental form because essentially normal connection is determined by the second fundamental form using the formula:
    $$\nabla^\perp_X\xi=\wt\nabla_X\xi+A_\xi(X)$$
    for $X\in \Gamma(TM)$ and $Y\in \Gamma(T^\perp M)$, where $A$ is the shape operator.  
\end{enumerate}
\end{rem}

We next come to prove Theorem \ref{thm-main}, the global existence of isometric existence with prescribed second fundamental form. 

\begin{proof}[Proof of Theorem \ref{thm-main}] For $x\in M$, let $\gamma_0,\gamma_1:[0,1]\to M$ be two smooth  curves joining $p$ to $x$. Since $M$ is simply connected, there is a smooth map $\Phi:[0,1]\times [0,1]\to M$ such that
\begin{equation*}
\left\{\begin{array}{ll}\Phi(0,t)=\gamma_0(t)&\mbox{for $t\in [0,1]$}\\
\Phi(1,t)=\gamma_1(t)&\mbox{for $t\in [0,1]$}\\
\Phi(u,0)=p&\mbox{for $u\in [0,1]$}\\
\Phi(u,1)=x&\mbox{for $u\in [0,1]$}.
\end{array}\right.
\end{equation*}
Let $\gamma_u(t)=\Phi(u,t)$, $\wt\Phi(u,t)=\wt\gamma_u(t)$ for any $u\in [0,1]$. 
Let $e_1,e_2,\cdots,e_n$ and $e_{n+1},e_{n+2},\cdots,e_{n+s}$ be orthonormal bases of $T_pM$ and $V_p$ respectively, and let $\wt e_A=\varphi(e_A)$ for $A=1,2,\cdots,n+r$. Suppose that
\begin{equation*}
v_{\gamma_u}(t)=v_a(u,t)e_a.
\end{equation*}
Then,
\begin{equation*}
\wt v_{\gamma_u}(t)=v_a(u,t)\wt e_a.
\end{equation*}
Let $E_A(u,t)=P_0^t(\gamma_u)(e_A)$ for $A=1,2,\cdots,n+s$, and 
\begin{equation*}
h(E_a,E_b)=h_{ab}^\alpha(u,t)E_\alpha.
\end{equation*}
Then, 
\begin{equation*}
\wt h_{\gamma_u}(\wt e_a,\wt e_b)=h_{ab}^\alpha \wt e_\alpha.
\end{equation*}
Let $\tilde E_A(u,t)=D_0^t(\tilde\gamma_u)(\tilde e_A)$ for $A=1,2,\cdots,n+s$. 
Suppose that
\begin{equation*}
\frac{\p\Phi}{\p u}=U_aE_a,\ \nabla_\pu E_a=X_{ab}E_b\mbox{ and }D_{\pu}E_\alpha=X_{\alpha\beta}E_\beta.
\end{equation*}
Then, by the same computation as in the proof Theorem \ref{thm-g-Jacobi}, we know that $U$ and $X$ satisfy the following system of ODEs:
\begin{equation}\label{eq-g-Jacobi}
\begin{split}
\left\{\begin{array}{l}U''_a=R_{cadb}v_cv_dU_b+\p_tv_bX_{ba}+\p_u\p_tv_a\\
X'_{ab}=R_{abdc}v_dU_c\\
X'_{\alpha\beta}=R^V_{\alpha\beta ab}v_aU_b\\
U_a(u,0)=X_{ab}(u,0)=X_{\alpha\beta}(u,0)=0\\
U'_a(u,0)=\p_uv_a(u,0),
\end{array}\right.
\end{split}
\end{equation}
where $R_{abcd}=R(E_a,E_b,E_c,E_d)$ and $R^V_{\alpha\beta ab}=R^V(E_\alpha,E_\beta,E_a,E_b)$.

Moreover, suppose that
\begin{equation*}
\frac{\p\wt\Phi}{\p u}=\wt U_AE_A, \mbox{ and }\wt \nabla_{\pu}\wt E_A=\wt X_{AB}\wt E_B.
\end{equation*}
Then, by the same argument as in the proof of Theorem \ref{thm-Cartan} using Theorem \ref{thm-g-Jacobi}, we have
\begin{equation}\label{eq-solution}
\left\{\begin{array}{l}\wt U_a=U_a\\
\wt U_\alpha=0\\
\wt X_{ab}=X_{ab}\\
\wt X_{\alpha\beta}=X_{\alpha\beta}\\
\wt X_{a\alpha}=h_{ab}^\alpha U_b.
\end{array}\right.
\end{equation}
By the \eqref{eq-solution}, we know that 
$$\wt U_a(u,1)=U_a(u,1)=0 \mbox{ and } \wt U_\alpha(u,1)=0.$$ 
So $$\frac{\p\wt\Phi}{\p u}(u,1)=0$$
 which implies that $\wt \gamma_u(1)=\gamma_0(1)$ for $u\in [0,1]$. Thus, $f$ is well defined. Then, by the same argument as in the proof of Theorem \ref{thm-Cartan}, we know that $f$ is an isometric immersion and $f_*(E_a)=\wt E_a$.

Moreover, by the definition of $\tau_\gamma$ and that $f_*(E_a)=\wt E_a$, we have 
$$\tau_{\gamma_u}|_{V_x}:V_x\to T_{f(x)}^\perp M,\ \left(\tau_{\gamma_u}|_{V_x}\right)(E_\alpha(u,1))=\wt E_\alpha(u,1). $$
Note that 
$$\p_uE_\alpha(u,1)=D_{\pu}E_\alpha(u,1)=X_{\alpha\beta}(u,1)E_\beta$$
and 
$$\p_u\wt E_\alpha(u,1)=\wt\nabla_{\pu}\wt E_\alpha(u,1)=\wt X_{\alpha\beta}(u,1)\wt E_\beta+\wt X_{\alpha a}(u,1)\wt E_a=X_{\alpha\beta}(u,1)\wt E_\beta$$
by that
$$\wt X_{a\alpha}(u,1)=-\wt X_{a\alpha}(u,1)=-h_{ab}^\alpha U_b(u,1)=0$$
and $X_{\alpha\beta}=\wt X_{\alpha\beta}$ using \eqref{eq-solution}. Therefore $\tau_{\gamma_u}|_{V_x}$ is independent of $u\in [0,1]$ and thus 
$\wt f$ is well-defined. The other properties $\wt f$ can be verified by the same argument as in the proof of Theorem \ref{thm-Cartan}. This completes the proof of the theorem.
\end{proof}

\end{document}